\documentclass[a4paper, BCOR=0.0mm, DIV=calc]{scrartcl}
\usepackage[T1]{fontenc}
\usepackage[utf8]{inputenc}
\usepackage[ngerman]{babel}
\usepackage[osf,sc]{mathpazo}
\usepackage{textcomp}
\usepackage{microtype}
\usepackage{amsmath, amssymb, amsfonts}
\usepackage{graphicx, color, wrapfig}
\usepackage{tabularx}
\KOMAoptions{twoside=false, twocolumn=false, headinclude=false, footinclude=false, mpinclude=false, pagesize=auto}
\recalctypearea

\setkomafont{section}{\mdseries\scshape\Large}
\setkomafont{title}{\mdseries\scshape}

\begin{document}
\title{Beweise bestimmter zahlentheoretischer Theoreme\footnote{
Originaltitel: "`Theorematum quorundam arithmeticorum demonstrationes "', erstmals publiziert in "`\textit{Commentarii academiae scientiarum Petropolitanae} 10, 1747, pp. 125 - 146"', Nachdruck in "`\textit{Opera Omnia}: Series 1, Volume 2, pp. 38 - 58"' und "`\textit{Commentat. arithm.} 1, 1849, pp. 24 - 34 [98a]"' , Eneström-Nummer E98, übersetzt von: Alexander Aycock, Textsatz: Artur Diener, im Rahmen des Projektes "`Euler-Kreis Mainz"' }}
\author{Leonhard Euler}
\date{}
\maketitle
Zahlentheoretische Theoreme, von welcher Art Fermat und andere viele entdeckt haben, sind umso größerer Aufmerksamkeit würdig, je geheimer die Gültigkeit derer und je schwerer sie zu beweisen sind. Fermat hat freilich eine hinreichend große Menge solcher Theoreme hinterlassen, hat aber niemals Beweise erörtert, auch wenn er vehement über deren Gültigkeit behauptete, dass sie sehr sicher feststeht. Besonders ist also zu bedauern, dass seine Schriften so sehr verloren sind, dass man alle Beweise noch immer nicht kennt. Ähnlich ist auch die Art der im Allgemeinen bekannten Propositionen, bei denen man behauptet, dass weder die Summe noch die Differenz zweier Biquadrate ein Quadrat ergeben kann, obwohl nämlich über deren Gültigkeit niemand zweifelt, existiert dennoch nirgendswo ein, wie weit es mir freilich bekannt ist, strenger Beweis, außer ein von Frénicle einst geschriebenes Büchlein, dessen Titel \textit{Traité des Traingles rectangles en nombres} ist. Der Autor beweist aber hier unter anderem, dass in keinem rechtwinkligen Dreieck, dessen Seiten durch rationale Zahlen ausgedrückt werden, die Fläche ein Quadrat sein kann, woher leicht die Gültigkeit der erwähnten Propositionen über die Summe und Differenzen zweier Biquadrate abgeleitet wird. Aber dieser Beweis ist so sehr in die Eigenschaften der Dreiecke eingehüllt, dass, wenn man nicht die größte Aufmerksamkeit aufbringt, es kaum klar eingesehen werden kann. Deswegen glaube ich, dass es der Mühe Wert sein wird, wenn ich die Beweise dieser Propositionen von rechtwinkligen Dreiecken losgemacht habe und sie analytisch und klar vorgelegt habe. Umso größere Nützlichkeit bringt aber dieses mein Unterfangen, je mehr andere um vieles schwerere Theoreme aus diesen gefunden werden können. Hierauf bezieht sich natürlich jenes gefeierte Theorem des Fermat, in dem er sagt, dass keine Dreieckszahl ein Biquadrat sein kann, außer der Einheit, den Beweis welches es mir aus jenen geglückt ist zu beweisen. Umso schwerer aber scheint der Beweis, weil die Proposition von einer Ausnahme abhängig ist und sich nur auf ganze Zahlen bezieht; bei gebrochenen Zahlen kann es nämlich auf unendlich viele Arten erreicht werden, dass $\frac{x(x+1)}{2}$ ein Biquadrat wird. Um also dieses und einige andere Theoreme zu zeigen, wird es nötig sein, gewisse Lemmata vorauszuschicken, auf die die folgenden Beweise gestützt sind; vor allem muss aber erinnert worden sein, dass für mich alle Buchstaben immer ganze Zahlen bezeichnen.
\section*{Lemma 1}
Das Produkt aus zwei oder mehreren zueinander primen Zahlen kann weder ein Quadrat noch ein Kubus noch irgendeine andere Potenz sein, wenn die einzelnen Faktoren nicht Quadrate oder Kuben oder andere Potenzen von solcher Art sind.

Der Beweis dieses Lemmas ist leicht und auch schon von Euklid angegeben worden, sodass es überflüssig wäre, diesen hier zu erörtern.
\section*{Lemma 2}
Wenn $a^2 + b^2$ ein Quadrat war und $a$ und $b$ zueinander prime Zahlen sind, wird $a = pp - qq$ und $b = 2pq$ sein, während $p$ und $q$ zueinander prime Zahlen sind, wobei die eine gerade, die andere ungerade ist.
\section*{Beweis}
Weil $a^2 + b^2$ ein Quadrat ist, setze man seine Wurzel gleich $a + \frac{bq}{p}$, wo ich den Bruch $\frac{a}{p}$ in kleinsten Termen ausgedrückt habe, sodass $p$ und $q$ zueinander prime Zahlen sind. Nachdem die Gleichung aufgestellt worden ist, wird $a^2 + b^2 = a^2 + \frac{2abq}{p} + \frac{bbqq}{pb}$ sein. Daher wird
\[
	a : b = (pp - qq) : 2pq.
\]
Die Zahlen $pp-qq$ und $2pq$ sind zueinander entweder prim oder haben den gemeinsamen Teiler $2$. In jenem Fall also, in dem $pp-qq$ und $2pq$ zueinander prime Zahlen sind, was passiert, wenn die eine der Zahlen $p$ und $q$ gerade war, die andere ungerade, ist es nötig, dass
\[
	a = pp - qq \quad \text{und} \quad b = 2pq
\]
ist, weil $a$ und $b$ als zueinander prime Zahlen festgesetzt werden. In dem Fall aber, in dem die Zahlen $pp-qq$ und $2pq$ den gemeinsamen Teiler $2$ haben, was sein wird, wenn jede der beiden Zahlen $p$ und $q$ ungerade war (eine der beiden kann nämlich nicht gerade sein, weil sie als zueinander prim gesetzt wurden), wird $a = \frac{pp-qq}{2}$ und $b=pq$ sein. Man setze aber $p+q = 2r$ und $p-q = 2s$; es werden $r$ und $s$ zueinander prime Zahlen sein und die eine derer ist gerade, die andere ungerade, woher
\[
	a = 2rs \quad \text{und} \quad b = rr - ss
\]
wird, welcher Ausdruck, weil er mit dem ersten übereinstimmt, anzeigt, wenn $aa+bb$ ein Quadrat war und $a$ und $b$ zueinander prime Zahlen sind, dass die eine derer die Differenz zweier zueinander primer Quadrate ist, von denen die eine gerade ist, die ungerade, die andere Zahl, aber gleich dem doppelten Produkt aus den Wurzeln dieser Quadrate ist. Das heißt, dass
\[
	a = pp - qq\quad \text{und} \quad b = 2pq
\]
ist, während $p$ und $q$ zueinander prime Zahlen sind, von denen die eine gerade ist, die andere ungerade ist.

\hfill\textsc{Q.E.D.}
\section*{Korollar 1}
Wenn also die Summe zweier zueinander primer Quadrate ein Quadrat war, ist es notwendig, dass das eine Quadrat gerade ist, dass andere aber ungerade; daraus folgt, dass die Summe zweier ungerader Quadrate kein Quadrat sein kann.
\section*{Korollar 2}
Wenn also $aa+bb$ ein Quadrat ist, wird die eine der Zahlen $a$ und $b$, z.\,B. $a$, ungerade sein, die andere $b$ aber gerade. Ungerade wird aber $a=pp-qq$ sein und gerade $b=2pq$.
\section*{Korollar 3}
Weil weiter die eine der Zahlen $p$ und $q$ gerade ist, die andere ungerade, wird $b$ eine mit $2$ multiplizierte Zahl oder durch $4$ teilbar sein. Wenn darauf weder $p$ noch $q$ durch $3$ teilbar ist, ist es notwendig, dass entweder $p-q$ oder $p+q$ eine Teilung durch $3$ zulässt. Daher folgt, dass die eine der Zahlen $a$ und $b$, die Summe welcher Quadrate ein Quadrat ergeben soll, durch $3$ teilbar sind.
\section*{Korollar 4}
Weil $a = pp-qq$ und $b=2pq$ ist, sieht man, wenn $aa+bb$ ein Quadrat festsetzt, leicht ein, dass die Zahlen $p$ und $q$ kleiner als $a$ und $b$ sind. Weil ja nämlich $a = (p+q)(p-q)$ ist, wird $a>p+q$ sein, wenn nicht $p-q=1$ ist; und wegen $b=2pq$ wird $p$ größer als $p$ oder $q$ sein. Auf diese Art werden also Zahlen $a$ und $b$ größer sein als $p$ und $q$. Es würde freilich $a=0$ werden, wenn $p=q$ wäre, aber dieser Fall hat keine Geltung, weil $p$ und $q$ als zueinander prime Zahlen festgesetzt worden sind und die eine derer gerade, die andere ungerade ist.
\section*{Bemerkung}
Beim Beweis dieses Lemmas folgt aus der Analogie $a:b = (pp-qq):2pq$ daher, dass $a=pp - qq$ und $b=2pq$ ist, weil $a$ und $b$ zueinander prime Zahlen sind und in gleicher Weise die Zahlen $pp-qq$ und $2pq$. Wenn nämlich $a:b = c:d$ war und die Zahlen $a$ und $b$ wie die Zahlen $c$ und $d$ zueinander prim sind, ist notwendig, dass $a=c$ und $b=d$ ist, wie leicht aus der Natur der Verhältnisse feststeht.
\section*{Lemma 3}
Wenn $aa-bb$ ein Quadrat war, während $a$ und $b$ zueinander prime Zahlen waren, wird $a = pp+qq$ sein und entweder $b=pp-qq$ und $b=2pq$, wo die Zahlen $p$ und $q$ zueinander prim sind und die eine derer gerade ist, die andere ungerade.
\section*{Beweis}
Weil $aa-bb$ ein Quadrat ist, setze man $a^2 - b^2 = c^2$ und es wird $s^2 = b^2 + c^2$ sein und die Zahlen $b$ und $c$ zueinander prime Zahlen. Weil also durch Korollar $1$ des vorhergehenden Lemmas die eine der Zahlen $b$ und $c$ gerade ist, die andere ungerade, ist notwendig, dass $a$ eine ungerade Zahl ist; $b$ wird in der Tat entweder gerade oder ungerade sein. Es sei zuerst $b$ ungerade und $c$ gerade, es wird durch das vorhergehende Lemma $b = pp-qq$ und $c = 2pq$ sein, während $p$ und $q$ zueinander prime Zahlen sind, von welchen die eine gerade ist, die andere ungerade. Daher aber wird $a=pp+qq$ sein. Aber wenn $b$ eine gerade Zahl war und $c$ ungerade, wird $b = 2pq$ und $c = pp-qq$ sein, woher erneut $a=pp+qq$ wird. Wenn deshalb $aa-bb$ ein Quadrat war, wird
\[
	a = pp+qq \quad \text{und sogar} \quad b = pp-qq \quad \text{oder} \quad b = 2pq
\] 
sein.

\hfill\textsc{Q.E.D.}
\section*{Korollar 1}
Wenn also die Differenz zweier Quadrate eine Quadratzahl ist, muss das größere Quadrat eine ungerade Zahl sein, wenn freilich jene Quadrate zueinander prime Zahlen waren.
\section*{Korollar 2}
Auf die gleiche Weise sieht man weiter ein, dass die Zahlen $p$ und $q$ kleiner als die Zahlen $a$ und $b$ sind, weil $a = pp+qq$ ist und $b$ entweder gleich $pp-qq$ oder gleich $2pq$ ist.
\section*{Korollar 3}
Wenn $aa-bb = cc$ war, ist eine der Zahlen $a$, $b$, $c$ immer durch $5$ teilbar. Denn weil $a = pp+qq$, $b=pp-qq$ und $c=2pq$ ist, ist entweder die eine der Zahlen $p$ und $q$ durch $5$ teilbar oder keine; in jenem Fall aber wird $c$ durch $5$ teilbar. In diesem Fall aber werden $pp$ und $qq$ Zahlen der Form solcher Art $5n\pm 1$ sein, also wird entweder $pp-qq$ oder $pp+qq$ durch $5$ teilbar sein.
\section*{Theorem 1}
Die Summe zweier Biquadrate wie $a^4 + b^4$ kann kein Quadrat sein, außer die eine der Zahlen verschwindet.
\section*{Beweis}
Beim Beweisen dieses Theorems werde ich so vorgehen, dass ich zeige, wie groß auch immer die Zahlen $a$ und $b$ waren, dass dann immer kleinere Zahlen anstelle von $a$ und $b$ angegeben werden können und man schließlich zu kleinsten Zahlen solche nicht gegeben sind, die Summe welcher Biquadrate ein Quadrat fesetsetzen würde, wird zu schlussfolgern sein, dass auch unter den größten Zahlen solche nicht existieren. Wir wollen also setzen, dass $a^4 + b^4$ ein Quadrat ist und $a$ und $b$ zueinander prime Zahlen sind; wenn sie nämlich nicht prim wären, könnte sie durch Teilung auf prime zurückgeführt werden. Es sei $a$ eine ungerade Zahl, $b$ aber eine gerade, weil notwendig die eine gerade, die andere ungerade sein muss. Es wird also
\[
	aa = pp-qq \quad \text{und} \quad bb=2pq
\]
sein und die Zahlen $p$ und $q$ werden zueinander prim sein und deren eine wird gerade, deren andere ungerade sein. Weil aber $aa = pp-qq$ ist, ist nötig, dass $p$ eine ungerade Zahl ist, weil andernfalls $pp-qq$ kein Quadrat sein könnte. Es wird also $p$ eine ungerade Zahl sein und $q$ eine gerade Zahl. Weil weiter $2pq$ ein Quadrat sein muss, ist notwendig, dass $p$ wie $2q$ ein Quadrat ist, weil $p$ und $2q$ zueinander prime Zahlen sind. Damit aber $pp-qq$ ein Quadrat ist, ist notwendig, dass
\[
	p = mm+nn \quad \text{und} \quad q = 2mn
\]
ist, während wiederum $m$ und $n$ zueinander prime Zahlen sind und die eine gerade, die andere ungerade ist. Weil weiter $2pq$ ein Quadrat sein muss, ist notwendig, dass so $p$ wie $2q$ ein Quadrat ist, weil $p$ und $2q$ zueinander prime Zahlen sind. Damit aber $pp-qq$ ein Quadrat ist, ist nötig dass
\[
	p = mm + nn \quad \text{und} \quad q = 2mn
\]
ein Quadrat ist, während wiederum $m$ und $n$ zueinander prime Zahlen sind und die eine derer gerade, die andere ungerade ist. Aber weil ja $2q$ ein Quadrat ist, wird $4mn$ oder $mn$ ein Quadrat sein; daher werden $m$ wie $n$ Quadrate sein. Nachdem also
\[
	m = xx\quad \text{und} \quad n = yy
\]
gesetzt wurde, wird
\[
	p = m^2 + n^2 = x^4 + y^4
\]
sein, was in gleicher Weise ein Quadrat sein müsste. Daher folgt also, wenn $a^4 + b^4$ ein Quadrat wäre, dass dann auch $x^4 + y^4$ ein Quadrat sein würde; es ist aber klar, dass die Zahlen $x$ und $y$ weit kleiner sein werden als $a$ und $b$. Auf dem gleichen Wege werden aus den Biquadraten $x^4 + y^4$ wieder kleinere Quadrate entstehen, die Summe welcher ein Quadrat wäre, und indem man so verfährt, wird man schließlich zu kleinsten Biquadraten in ganzen Zahlen gelangen. Weil also kleinste Biquadrate nicht gegeben sind, deren Summe ein Quadrat erzeugen würde, ist offensichtlich, dass auch in größten Zahlen solche nicht gegeben sind. Wenn aber bei einem geraden der Biquadrate das eine gleich $0$ ist, wird bei allen anderen übrigen geraden das andere verschwinden, sodass daher keine neuen Fälle entstehen.

\hfill\textsc{Q.E.D.}
\section*{Korollar 1}
Weil also die Summe zweier Biquadrate kein Quadrat sein könnte, können um vieles weniger zwei Biquadrate zusammen ein Biquadrat erzeugen.
\section*{Korollar 2}
Obwohl dieser Beweis sich nur auf ganze Zahlen bezieht, erledigt man dennoch auch durch ihn, dass nicht einmal in gebrochenen Zahlen zwei Biquadrate beschafft werden können, deren Summe ein Quadrat wäre. Denn wenn $\frac{a^4}{m^4} + \frac{b^4}{n^4}$ ein Quadrat wäre, dann wäre in ganzen Zahlen $a^4 n^4 + b^4 m^4$ ein Quadrat, was durch den Beweis selbst nicht passieren kann.
\section*{Korollar 3}
Aus demselben Beweis lässt sich berechnen, dass keine Zahlen $p$ und $q$ solcher Art gegeben sind, dass $p$, $2q$ und $pp-qq$ Quadrate sind; wenn solche nämlich existierten, dann würde man Werte für $a$ und $b$ haben, die $a^4 + b^4$ zu einem Quadrat machen; es wäre nämlich $a = \sqrt{pp-qq}$ und $b = \sqrt{2pq}$.
\section*{Korollar 4}
Nachdem also $p=xx$ und $2q=4yy$ gesetzt wurden, wird $pp-qq = x^4 - 4y^4$ sein. Es kann also insgesamt nicht geschehen, dass $x^4 - 4y^4$ ein Quadrat ist. Und es wird auch $4x^4 - y^4$ kein Quadrat sein können; es wäre nämlich $16x^4 - 4y^4$ ein Quadrat, welcher Fall wegen des Biquadrates $16x^4$ auf den ersten zurückfällt.
\section*{Korollar 5}
Es folgt daher auch, dass $ab(a^2 + b^2)$ niemals ein Quadrat sein kann. Wegen der zueinander primen Faktoren $a$, $b$, $a^2+b^2$ müssten die einzelnen nämlich Quadrate sein, was nicht passieren kann.
\section*{Korollar 6}
Ähnlich werden auch solche zueinander primen Zahlen $a$ und $b$ nicht gegeben sein, die das Quadrat $2ab(aa-bb)$ erzeugen würden. Dies folgt aus Korollar $3$, wo bewiesen wurde, dass keine Zahlen $p$ und $q$ gegeben sind, sodass $p$, $2q$, $pp-qq$ Quadrate wären. All dies gilt aber auch für zueinander nicht prime Zahlen und sogar für gebrochene Zahlen durch Korollar $2$.
\section*{Theorem 2}
Die Differenz zweier Biquadrate wie $a^4 - b^4$ kann kein Quadrat sein, wenn nicht entweder $b=0$ oder $b=a$ war.
\section*{Beweis}
Ich werde dieses Theorem auf die gleiche Weise beweisen wie das vorhergehende. Es seien also die Biquadrate schon auf kleinste Terme reduziert und wir wollen setzen, dass $a^4 - b^4$ ein Quadrat ist; es wird aber $a$ eine ungerade Zahl sein, $b$ aber entweder gerade oder ungerade sein.
\section*{Fall 1}
Es sei zunächst $b$ eine gerade Zahl; es wird
\[
	a^2 = pp + qq \quad \text{und} \quad b^2 = 2pq
\]
sein, während $p$ und $q$ zueinander prime Zahlen sind und die eine $p$ gerade, die andere $q$ ungerade ist. Wegen $b^2 = 2pq$ werden also $2p$ und $q$ Quadrate sein müssen. Weil weiter $pp+qq$ gleich $a^2$ wird, wird
\[
	q = mm - nn \quad \text{und} \quad p=2mn
\]
sein, während $m$ und $n$ zueinander prime Zahlen sind. Weil aber $2p$ ein Quadrat ist, wird dieses $4mn$ sein, das heißt, $mn$ ist ein Quadrat; und daher sind $m$ und $n$ einzeln Quadrate. Nachdem also
\[
	m = x^2 \quad \text{und} \quad n=y^2
\]
gesetzt wurde, wird
\[
	q = x^4 - y^4
\]
sein; weil dort die eine der Zahlen $m$ und $n$ gerade ist, die andere ungerade, wird auch die eine der Zahlen $x$ und $y$ gerade sein, die andere ungerade. Aber wegen des Quadrates $q$ wird $x^4 - y^4$ ein Quadrat sein, wo $x$ eine ungerade Zahl sein wird, $y$ aber eine gerade. Wenn deshalb $a^4 - b^4$ ein Quadrat war, wird auch $x^4 - y^4$ ein Quadrat sein, während $x$ und $y$ weit kleiner als $a$ und $b$ sind. Weil also in kleinsten Zahlen keine zwei Biquadrate gegeben sind, die als Differenz ein Quadrat haben, werden sie auch in größten nicht gegeben sein, zumindest in dem Fall, in dem das kleinere Biquadrat eine gerade Zahl war. 

\hfill\textsc{Q.E. Unum}
\section*{Fall 2}
Es sei nun $b$ eine ungerade Zahl und es wird
\[
	a^2 = pp+qq \quad \text{und} \quad b^2 = pp-qq
\]
sein, während $p$ und $q$ zueinander prime Zahlen sind und die eine von diesen gerade, die andere ungerade ist. Weil aber $pp-qq$ ein Quadrat ist, wird $p$ eine ungerade Zahl sein und deshalb $q$ eine gerade. Nachdem aber $a^2$ und $b^2$ miteinander multipliziert wurden, wird $a^2 b^2 = p^4 - q^4$ hervorgehen, welcher Ausdruck durch den ersten Fall kein Quadrat sein kann und daher $a^2 b^2$ nicht gleich werden kann. Die Differenz zweier Biquadrate kann also auf keine Weise ein Quadrat sein, wenn nicht entweder beide gleich sind oder eines gleich $0$. 

\hfill\textsc{Q.E. Alterum Dem.}
\section*{Korollar 1}
Weil $a^2 = pp+qq$ und $b^2 = 2pq$ ist und so $q = mm-nn$ und $p=2mn$ und weiter $m=x^2$ und $n=y^2$, wird $a^2 = (x^4 + y^4)^2$ und $b^2 = 4x^2 y^2 (x^4 - y^4)$ sein. Daraus wird man $a = x^4 + y^4$ und $b = 2xy\sqrt{x^4 + y^4}$ haben.
\section*{Korollar 2}
Wenn also in kleinsten Zahlen $x$ und $y$ solche gegeben wären, die Differenz welcher Biquadrate ein Quadrat festsetzen würden, dann könnten aus denen sofort um vieles größere Zahlen, die sich der selben Eigenschaft erfreuen, $a$ und $b$ gefunden werden.
\section*{Korollar 3}
Daher erkennt man klarer, dass der Fall, in dem zwei Biquadrate entweder gleich sind oder das eine gleich $0$ ist, keine neuen Fälle liefert; nachdem nämlich entweder $x=y$ und $y=0$ gesetzt worden ist, wird zugleich $b=0$ werden, woher man die Weite des Beweises umso mehr erfasst.
\section*{Korollar 4}
Aus dem Beweis folgt weiter, dass keine Zahlen $p$ und $q$ der Formen gegeben sind, dass $2p$, $q$ und $qq+pp$ Quadrate sind. Nachdem also $2p=4xx$ und $q=yy$ gesetzt wurde, wird kein Quadrat in dieser Form $4x^4 + y^4$ sein können.
\section*{Korollar 5}
Aus diesen Formeln folgt auch, dass weder $ab(aa-bb)$ noch $2ab(aa+bb)$ jeweils Quadrate sein können, was nicht nur gilt, wenn $a$ und $b$ zueinander prime Zahlen sind, sondern auch wenn sie zusammengesetzt und sogar gebrochen waren. Denn Brüche von solcher Art werden leicht ganze Zahlen und ganze Zahlen auf zueinander prime Zahlen zurückgeführt.
\section*{Korollar 6}
In diesen zwei Propositionen ist gezeigt worden, dass die folgenden neun Ausdrücke niemals Quadrate sein können:
\[
\begin{array}{llll}
\text{I.} & a^4 + b^4 &\qquad \text{VI.} & a^4 - b^4 \\
\text{II.} & a^4 - 4b^4 &\qquad \text{VII.} & 4a^4 + b^4 \\
\text{III.} & 4a^4 - b^4 &\qquad \text{VII.} & ab(aa-bb) \\
\text{IV.} & ab(aa+bb) &\qquad \text{IX.} & 2ab(aa+bb) \\
\text{V.} & 2ab(aa-bb) &\qquad \text{X.} & 2a^4 \pm 2b^4 
\end{array}
\] 
Ich habe den $10.$ Ausdruck daher hinzugefügt, weil seine Gültigkeit bald gezeigt werden wird.
\section*{Theorem 3}
Die Summe zweier Biquadrate, die $2$-mal genommen wurden, wie $2a^4 + 2b^4$, kann kein Quadrat sein, wenn nicht $a=b$ war.
\section*{Beweis}
Ich setze zuerst, dass $a$ und $b$ zueinander prime Zahlen sind; denn wenn sie solche nicht wären, könnte die Formel durch Teilung darauf zurückgeführt werden. Man erkennt aber leicht, dass jede der beiden Zahlen $a$ und $b$ ungerade sein muss; wenn nämlich eine gerade wäre, dann würde $2a^4 + 2b^4$ eine mit $2$ multiplizierte ungerade Zahl werden, die kein Quadrat sein kann. Weiter stimmt diese Form mit dieser $(aa+bb)^2 + (aa-bb)^2$ überein, welche daher gezeigt werden muss, kein Quadrat sein zu können, wenn nicht $a=b$ war. Aber wegen der ungeraden Zahlen $a$ und $b$ werden $a^2 + b^2$ und $a^2 - b^2$ gerade Zahlen sein, jene natürlich ungerade, diese aber gerade. Man ist also auf diese Form $\left( frac{aa+bb}{2}\right)^2 + \left(\frac{aa-bb}{2}\right)^2$ geführt worden, in welcher $\frac{aa+bb}{2}$ und $\frac{aa-bb}{2}$ zueinander prime Zahlen sind, jene ungerade, diese aber gerade; wenn deshalb die vorgelegte Form ein Quadrat wäre, würde
\[
	\frac{aa+bb}{2} = pp-qq\quad \text{und} \quad \frac{aa-bb}{2} = 2pq
\]
sein, woher man $a^2 = pp+2pq + qq$ und $b^2 = pp-2pq-qq$ findet, die Differenz welcher Ausdrücke
\[
	4pq = aa-bb
\]
ist; und daher wird $a+b = \frac{2mp}{n}$ und $a-b = \frac{2nq}{m}$ sein, woher
\[
	a = \frac{mp}{n} + \frac{nq}{m} \quad \text{und} \quad b = \frac{mp}{n} - \frac{nq}{m}
\]
würde. Nachdem aber diese Substitution gemacht worden ist, wird
\[
	\frac{mm}{nn}pp + \frac{nn}{mm}qq = pp-qq \quad\text{und}\quad \frac{pp}{qq} = \frac{nn(mm+nn)}{mm(nn-mm)} - \frac{nn(n^4 - m^4)}{mm(nn-mm)^2}
\]
sein. Es würde also $n^4 - m^4$ ein Quadrat sein müssen, was durch das vorhergehende Theorem nicht geschehen werden kann.

\hfill\textsc{Q.E.D.}
\section*{Korollar 1}
Wenn also $a$ und $b$ ungerade Zahlen waren, kann auch $2ab(aa+bb)$ kein Quadrat sein; es würden nämlich $a$, $b$ und $2aa + 2bb$ Quadrate sein müssen, was durch das Theorem nicht geschehen kann.
\section*{Korollar 2}
Der Beweis hätte also auch aus der neunten Formel $2ab(aa+bb)$ gebildet werden können; aber dort war die eine der Zahlen $a$ und $b$ gerade gesetzt worden, die andere ungerade; Wenn das auch nichts behindert, nützte es dennoch, einen gesonderten Beweis gegeben zu haben.
\section*{Korollar 3}
Durch diesen Beweis wird also die Gültigkeit der neuen Formeln mehr bestätigt, weil daher schon feststeht, dass $2ab(aa+bb)$ kein Quadrat sein kann, auch wenn die Zahlen $a$ und $b$ beide ungerade sind.
\section*{Korollar 4}
Kürzer aber kann die Gültigkeit dieses Theorems aus der Form $(a^2 + b^2)^2 + (a^2 - b^2)^2$ gezeigt werden, die daher kein Quadrat sein kann, weil $(a^2 + b^2)^2 - (a^2 - b^2)^2$ ein Quadrat ist. Es kann aber nicht geschehen, dass die Summe zweier Quadrate ein Quadrat ist, wenn die Differenz derselben Quadrate ein Quadrat war. Wenn nämlich so $pp+qq$ wie $pp-qq$ ein Quadrat wäre, wäre $p^4 - q^4$ ein Quadrat, was nicht sein kann.
\section*{Korollar 5}
Auf die gleiche Weise kann $a^4 - 6aabb + b^4$ kein Quadrat sein. Es ist nämlich
\[
	a^4 - 6aabb + b^4 = (aa-bb)^2 - 4aabb,
\]
was die Differenz zweier Quadrate solcher Art ist, deren Summe ein Quadrat gibt.
\section*{Korollar 6}
Und auf die gleiche Weise kann $a^4 + 6a^2 b^2 + b^4$ kein Quadrat sein, weil es gleich $(a^2 + b^2)^2 + 4aabb$ ist, die Summe welcher Quadrate kein Quadrat sein kann, weil die Differenz derselben $(a^2 + b^2)^2 - 4aabb$ ein Quadrat ist.
\section*{Theorem 4}
Das Doppelte der Differenz zweier Biquadrate wie $2a^4 - 2b^4$ kann kein Quadrat sein, wenn nicht $a=b$ war.
\section*{Beweis}
Wir wollen setzen, dass die Zahlen $a$ und $b$ zueinander prime Zahlen sind und $2a^4 - 2b^4$ ein Quadrat ist; es werden $a$ und $b$ ungerade Zahlen sein. Es wäre als $2(a-b)(a+b)(aa+bb)$ ein Quadrat und daher auch sein $16.$ Teil $\left(\frac{a-b}{2}\right)\left(\frac{a+b}{2}\right)\left(\frac{aa+bb}{2}\right)$; weil diese Faktoren zueinander prim sind, würden die einzelnen Quadrate  sein müssen. Es sei also
\[
	\frac{a-b}{2} = pp \quad \text{und} \quad \frac{a+b}{2} = qq;
\]
es wird
\[
	a = pp+qq\quad \text{und} \quad b = qq-pp
\]
sein, woher
\[
	\frac{aa+bb}{2} = p^4 + q^4
\]
wird. Weil also $p^4 + q^4$ kein Quadrat sein kann, kann auch $\frac{aa+bb}{2}$ und daher $2a^4 - 2b^4$ kein Quadrat sein.

\hfill\textsc{Q.E.D.}
\section*{Theorem 5}
Weder $ma^4 - m^3 b^4$ noch $2ma^4 - 2m^3 b^4$ kann ein Quadrat sein.
\section*{Beweis}
Wir wollen setzen, dass $a$ und $b$ zueinander prime Zahlen sind und $m$ eine weder quadratische noch eine durch ein Quadrat teilbare Zahl ist; wenn nämlich $m$ durch ein Quadrat teilbar wäre, dann könnte der quadratische Faktor durch Teilung weggeschafft werden. Man setze weiter, dass $m$ eine zu $a$ wie zu $b$ prime Zahl ist; es werden wegen
\[
	ma^4 - m^3 b^4 = m(aa-mbb)(aa+mbb)
\]
die ganzen Faktoren zueinander prim sein und daher müssten die einzelnen Quadrate sein. Nachdem also $m=pp$ gesetzt wurde, müsste $(aa-ppbb)(aa+ppbb)$ ein Quadrat sein, was nicht geschehen kann.

Auf die gleiche Weise wird wegen $2ma^4 - 2m^3 b^4 = 2m(aa-mbb)(aa+mbb)$ und der entweder zueinander primen Faktoren oder die die zwei als gemeinsamen Faktor haben entweder $2m$ oder $m$ ein Quadrat sein; im ersten Fall müsste aber für $2m = 4pp$ gesetzt $a^4 - 4p^4 b^4$ ein Quadrat sein, was in gleicher Weise nicht passieren kann. Wenn aber $m=pp$ war, dann würde $2a^4 - 2p^4 b^4$ ein Quadrat sein, was durch das vorhergehende Theorem nicht passieren kann.

Aber wenn $m$ in Bezug auf $a$ nicht prim war, wollen wir $m=rs$ und $a=rc$ setzen, wo zu bemerken ist, dass $r$ und $s$ zueinander prime Zahlen sind, weil $m$ festgesetzt worden ist, keinen quadratischen Faktor zu haben. Die Quadrate müssten also diese Formen $r^5 sc^4 - r^3 s^3 b^4$ und $2r^5 sc^4 - 2r^3 s^3 b^4$ oder $r^3 sc^4 - rs^3 b^4$ und $2r^3sc^4 - 2r^3sb^4$ sein. Wegen der Faktoren dieser Formeln, die zueinander prim sind, müssten entweder $rs$ oder $2rs$ Quadrate sein und daher $r$ und $s$ oder $2s$ getrennt, woher Formeln entstehen würden, welche gezeigt worden sind, keine Quadrate sein zu können.

\hfill\textsc{Q.E.D.}
\section*{Korollar 1}
Formen dieser Art $mn(m^2 a^4 - n^2 b^4)$ und $2mn(m^2 a^4 - n^2 b^4)$ können also keine Quadrate sein, welche Zahlen auch immer anstelle von $m$, $n$, $a$ und $b$ angenommen werden.
\section*{Korollar 2}
Wenn also $maa+nbb$ ein Quadrat war, werden weder $m^2 naa - mn^2 bb$ noch $2m^2 naa - 2mn^2 bb$ Quadrate sein können. Und wenn $maa-nbb$ ein Quadrat war, werden weder $m^2 naa + mn^2 bb$ noch $2m^2 naa + 2mn^2 bb$ Quadrate sein können.
\section*{Korollar 3}
Wir wollen $maa + nbb = cc$ setzen; es wird $m = \frac{cc-nbb}{aa}$ sein; es wird also weder $n(cc-nbb)(cc-2nbb)$ noch $2n(cc-nbb)(cc-2nbb)$ ein Quadrat sein können. Und wenn $m=\frac{cc+nbb}{aa}$ war, dann wird auch keine dieser Formeln $n(cc+nbb)(cc+2nbb)$ und $2n(cc+nbb)(cc+2nbb)$ ein Quadrat sein können.
\section*{Korollar 4}
Wenn $c = \pm pp + nqq$ und $b = 2pq$ gesetzt wird, wird man die folgenden Formeln $n(p^4 \pm 6nppqq + n^2 q^4)$ und $2n(p^4 \pm 6nppqq + n^2 q^4)$ erhalten, die auf keine Weise ein Quadrat erzeugen können.
\section*{Theorem 6}
Weder $ma^4 + m^3 b^4$ noch $2m^4 + 2m^3 b^4$ kann ein Quadrat sein.
\section*{Beweis}
Ich sage zuerst, wenn $mp^2 - mq^2$ ein Quadrat war, dass dann weder $mp^2 + mq^2$ noch $2mp^2 + 2mq^2$ auf irgendeine Art ein Quadrat sein können; es würde nämlich entweder $m^2(p^4-q^4)$ oder $2m^2(p^4-q^4)$ entgegen dem schon Bewiesenen ein Quadrat sein. Wir wollen aber $m^2 p^2 - mq^2$ zum Quadrat machen, indem wir seine Wurzel gleich $\frac{(p-q)a}{b}$ setzen; es wird $mp + mq = \frac{a^2 p - a^2 q}{bb}$ sein, woher man $q = \frac{p(aa-mbb)}{aa+mbb}$ findet. Es sei also $p = a^2 + mb^2$; es wird $q = a^2 - mb^2$ sein und daher $p^2 + q^2 = 2a^4 + 2m^2 b^4$. Es wird also zuerst $mp^2 + mq^2 = 2ma^4 + 2m^3 b^4$ kein Quadrat sein können, weiter $2mp^2 + 2mq^2 = 4ma^4 + 4m^3 b^4$. Daraus berechnet man, dass weder $ma^4 + m^3 b^4$ noch $2ma^4 + 2m^3 b^4$ ein Quadrat sein kann.

\hfill\textsc{Q.E.D.}
\section*{Korollar}
In diesen zwei Theoremen ist also gezeigt worden, dass keine Zahlen in diesen Formen $ma^4 \pm m^3 b^4$ und $2ma^4 \pm 2m^3 b^4$ Quadrate sein können. In diesen Formeln aber sind alle vorhergehenden enthalten.
\section*{Theorem 7 \\ (des Fermat)}
Keine Dreickszahl in ganzen Zahlen kann ein Biquadrat sein außer der Einheit.
\section*{Beweis}
Alle Dreieckszahlen sind in dieser Form $\frac{x(x+1)}{2}$ enthalten. Es ist also zu beweisen, dass diese Form $\frac{x(x+1)}{2}$ niemals ein Biquadrat sein kann, wenn freilich anstelle von $x$ ganze Zahlen (ausgenommen $x=1$) eingesetzt werden. Es ist aber zu bemerken, dass entweder $x$ eine gerade Zahl ist oder eine ungerade; im ersten Fall muss also $\frac{x}{2}(x+1)$, im zweiten aber $x\frac{x+1}{2}$ ein Biquadrat sein; in jeden der beiden dieser Faktoren würden zwei Faktoren zueinander prim sein müssen und daher müsste jeder der beiden ein Biquadrat sein. Es sei also im ersten Fall $\frac{x}{2} = m^4$ oder $x = 2m^4$ und es wird $x+1 = 2m^4 + 1$ ein Biquadrat sein müssen. Im zweiten Fall aber sei $\frac{x+1}{2} = m^4$, dass $x = 2m^4 - 1$ ist, was genauso ein Biquadrat sein muss. Deswegen müsste $2m^4 \pm 1$ ein Biquadrat sein. Man setze also $2m^4 \pm 1 = n^4$; es wird $4m^4 = 2n^4 \mp 2$ sein; es muss also $2n^4 \pm 2$ gleich $4m^4$ sein, das heißt ein Quadrat. Oben wurde aber gezeigt, dass $2a^4 \pm 2b^4$ und daher auch $2n^4 \pm 2$ niemals als Quadrat sein kann außer im Fall $n=1$. Nachdem aber $n=1$ gesetzt wurde, wird $m$ entweder $0$ oder $1$ und $x$ entweder $0$ oder $1$. Es ist also keine ganze Zahl gegeben, die anstelle von $x$ eingesetzt $\frac{x(x+1)}{2}$ ein Biquadrat ergäbe, außer den Fällen $x=0$ und $x=1$. Deshalb existiert in ganzen Zahlen keine Dreieckszahl, die ein Biquadrat wäre, außer der Einheit und der $0$.

\hfill\textsc{Q.E.D.}
\section*{Korollar 1}
Wenn man $\frac{xx+x}{2} = y^4$ setzt, wird $4xx + 4x + 1 = 8y^4 + 1 = (2x + 1)^2$ sein. Daraus folgt, indem man ganze Zahlen anstelle von $y$ einsetzt,dass die Form $8y^4 + 1$ niemals ein Quadrat sein kann außer in den Fällen $y=0$ und $y=1$.
\section*{Korollar 2}
Wenn $8y^4 + 1 = z^2$ gesetzt wird, wird $16y^4 = 2z^2 - 2$ werden. Deshalb kann $2z^2 - 2$ niemals ein Biquadrat werden, welche ganze Zahl auch immer anstelle von $z$ eingesetzt wird, außer in den Fällen $z=1$ und $z=3$.
\section*{Theorem 8}
Die Summe dreier Biquadrate, von denen zwei einander gleich sind, oder eine Form solcher Art $a^4 + 2b^4$, kann kein Quadrat sein, wenn nicht $b=0$ war.
\section*{Beweis}
Wir wollen setzen, dass $a^4 + 2b^4$ ein Quadrat ist und seine Wurzel $a^2 + \frac{m}{n}b^2$ ist, wo $a$ und $b$ wie $m$ und $n$ zueinander prime Zahlen sein werden. Nachdem aber die Gleichung aufgestellt worden ist, wird $2n^2 b^2 = 2mna^2 + m^2 b^2$ sein und
\[
	\frac{b^2}{a^2} = \frac{2mn}{2n^2 - m^2},
\]
welcher Bruch entweder schon die einfachste Form hat oder durch eine Teilung durch $2$ auf die einfachste zurückführbar sein wird. Wir wollen also zuerst setzen, dass $2mn$ und $2n^2 - m^2$ zueinander prime Zahlen sind, was passiert, wenn $m$ eine ungerade Zahl ist, und es wird
\[
	b^2 = 2mn \quad \text{und}\quad a^2 = 2n^2 - m^2
\]
sein. Hier sind zwei Fälle zu entwickeln, von denen der eine der ist, wenn $n$ eine ungerade Zahl ist, der andere, wenn $n$ gerade ist; in jenem Fall, in dem $n$ ungerade ist, ist klar, dass wegen des ungeraden $m$ auch $2mn$ kein Quadrat sein kann, in diesem Fall aber, in dem $n$ eine gerade Zahl ist, kann $a^2 = 2n^2 - m^2$ oder $a^2 + m^2 = 2n^2$ wegen der ungeraden Zahlen $a$ und $m$ und $2n^2$ keine doppelt gerade Zahl werden. Es mögen also $2mn$ und $2n^2 - m^2$ den gemeinsamen Teiler $2$ haben, was passiert, wenn $m$ eine gerade Zahl ist, z.\,B. $m=2k$, und es wird $n$ eine ungerade Zahl sein; man wird also $\frac{b^2}{a^2} = \frac{4kn}{2n^2 - 4k^2}$ haben, wo $2kn$ und $n^2 - 2k^2$ zueinander prime Zahlen sein werden. Daher wird also wegen der genauso zueinander primen Zahlen $b^2$ und $a^2$ 
\[
	b^2 = 2kn \quad \text{und} \quad a^2 = n^2- 2k^2
\]
sein. Aber hier kann $2kn$ kein Quadrat werden, wenn nicht $k$ eine gerade Zahl war. Es sei also $k$ eine gerade Zahl und es werden $n$ wie $2k$ Quadrate sein müssen; es werde also $n=cc$ und $2k=4dd$, wo $c$ eine ungerade Zahl sein wird, und dadurch wird man
\[
	a^2 = c^4 - 8d^4
\]
haben. Um also zu finden, ob $c^4-8d^4$ ein Quadrat sein kann, wollen wir setzen, dass seine Wurzel $c^2 - \frac{2p}{q}d^2$ ist und es wird $2q^2 d^2 = pqc^2 - p^2 d^2$ sein oder
\[
	\frac{dd}{cc} = \frac{pq}{pp+2qq},
\]
wo wiederum $c$ und $d$ wie $p$ und $q$ zueinander prime Zahlen sind. Hier sind erneut zwei Fälle zu bemerken, ob $p$ eine gerade oder eine ungerade Zahl ist. Es sei also zuerst $p$ eine ungerade Zahl; man wird wegen der zueinander primen Zahlen $pq$ und $pp+2qq$
\[
	dd = pq \quad \text{und} \quad cc = pp+2qq
\]
haben. Es ist also notwendig, dass $p$ wie $q$ ein Quadrat ist; deshalb setze ich $p = x^2$ und $q = y^2$ und es wird
\[
	cc = x^4 + 2y^4
\]
sein; wenn daher $a^4 + 2b^4$ ein Quadrat wäre, dann wäre auch $x^4 + 2y^4$ ein Quadrat und die Zahlen $x$ und $y$ werden um einiges kleiner sein als $a$ und $b$; daraus könnten erneut kleinere gefunden werden, was in ganzen Zahlen nicht geschehen kann. Für den $2.$ Fall, in dem $p$ eine gerade Zahl ist, wollen wir $p=2r$ setzen und es wird $\frac{dd}{cc} = \frac{2qr}{2rr+2qq} = \frac{qr}{2rr+qq}$ und wegen des ungeraden $p$ werden $qr$ und $2rr+qq$ zueinander prime Zahlen sein. Es wird also
\[
	dd = qr \quad\text{und}\quad cc = 2rr+qq
\]
sein; daher muss jeder der beiden Zahlen $p$ und $q$ ein Quadrat sein; nachdem deshalb $q = x^2$ und $r=y^2$ gesetzt worden ist, wird
\[
	cc = 2y^4 + x^4
\]
werden; daher ist klar, wenn $a^4 + 2b^4$ ein Quadrat wäre, dass dann auch bei weit kleineren Zahlen die gleiche Form $x^4 + 2y^4$ ein Quadrat sein würde. Deshalb kann $a^4 + 2b^4$ kein Quadrat sein, wenn nicht $b=0$ war.

\hfill\textsc{Q.E.D.}
\section*{Korollar 1}
Weil wir ja $\frac{b^2}{a^2} = \frac{2mn}{2n^2 - m^2}$, nachdem $a^4 + 2b^4$ als ein Quadrat festgesetzt wurde, gefunden haben, folgt, dass $2mn(2n^2 - m^2)$ kein Quadrat sein kann, welche Zahlen auch immer anstelle von $m$ und $n$ eingesetzt werden.
\section*{Korollar 2}
Nachdem also $m=x^2$ und $n=y^2$ gesetzt wurde, wird diese Form $4y^4 - 2x^4$ kein Quadrat sein. Auf gleiche Weise wird für $2m = 4x^2$ und $n=y^2$ diese Form $2y^4 - 4x^4$ kein Quadrat sein. Und für $m=x^2$ und $2n=4y^2$ gesetzt kann diese Formel $8y^4 -x^4$ kein Quadrat sein.
\section*{Korollar 3}
Wenn allgemein $m=\alpha x^2$ und $n=\beta y^2$ wird, wird diese Formel $2\alpha\beta (2\beta^2 y^4 - \alpha^2 x^4)$ oder $4\alpha \beta^3 y^4 - 2\alpha^3\beta x^4$ hervorgehen, die auf keine Weise ein Quadrat sein können wird.
\section*{Theorem 9}
Wenn diese Form $a^4 + kb^4$ kein Quadrat sein kann, dann wird auch diese Form $2k\alpha \beta^3 y^4 - 2\alpha^3 \beta x^4$ auf keine Weise ein Quadrat ergeben können.
\section*{Beweis}
Wir wollen also setzen, dass die vorgelegte Form $a^4 + kb^4$ ein Quadrat ist und ihre Wurzel gleich $a^2 + \frac{m}{n}b^2$ ist; es wird $kn^2 b^2 = 2mna^2 + m^2 b^2$ sein und $\frac{b^2}{a^2} = \frac{2mn}{kn^2 - m^2}$. Weil also $a^4 + bk^4$ kein Quadrat sein kann, wird dann auch $\frac{2mn}{kn^2 - m^2}$ oder $2mn(kn^2 - m^2)$ kein Quadrat sein können. Es werde $m=\alpha x^2$ und $n=\beta y^2$; es wird $2\alpha \beta (k\beta^2 y^4 - \alpha^2 x^4)$ oder $2k\alpha \beta ^3 y^4$ hervorgehen, welche Formel deshalb kein Quadrat sein kann, welche positiven oder negativen Zahlen auch immer anstelle von $\alpha$ und $\beta$ eingesetzt werden.

\hfill\textsc{Q.E.D.}
\section*{Korollar 1}
Es werde nun entweder $\alpha$ oder $\beta$ negativ, dass diese Form $2\alpha^3 \beta x^4 - 2k\alpha \beta^3 y^4$ hervorgeht, und man setze $2\alpha^3 \beta = p^2$; es wird $\beta = \frac{p^2}{2\alpha^3}$ sein, woher jene Form in diese $p^2 x^4 - \frac{kp^6}{4\alpha^8 y^4}$ übergeht. Diese Formel $x^4 - 4ky^4$ kann also kein Quadrat sein, nachdem $4y^4$ für $\frac{p^4}{4\alpha^8 y^4}$ gesetzt wurde. Aus dieser Formel folgt weiter, dass dieser Ausdruck $2\alpha^3 \beta x^4 + 8k\alpha \beta^3 y^4$ kein Quadrat sein kann.
\section*{Korollar 2}
Man setze $2k\alpha\beta^3 = pp$ in der gefundenen Formel $2k\alpha \beta^3 y^4 - 2\alpha^3 \beta x^4$, dass $\alpha = \frac{pp}{2k\beta^3}$ ist; jene wird in diese $p^2 y^4 - \frac{p^6}{4k^3 p^8}x^4$ übergehen, aus welcher folgt, dass $a^4 - 4kb^4$ kein Quadrat sein kann; daher wird wie zuvor $2\alpha^3 \beta x^4 + 8k\alpha \beta^3 y^4$ kein Quadrat sein können.
\section*{Korollar 3}
Wenn also $a^4 + kb^4$ kein Quadrat sein kann, dann wird weder diese Formel
\[
	2k\alpha \beta^3 y^4 - 2\alpha^3 \beta x^4
\]
noch diese
\[
	\alpha^3 \beta x^4 + k\alpha \beta^3 y^4
\]
ein Quadrat sein, welche letztere aus den vorhergehenden Korollaren folgt, indem man $2\alpha$ anstelle von $\alpha$ schreibt.
\section*{Korollar 4}
Weil also $a^4 + b^4$ kein Quadrat sein kann, werden die beiden folgenden Formeln $\alpha^3 \beta x^4 + \alpha \beta^3 y^4$ und $2\alpha \beta^3 y^4 - 2\alpha^3 \beta x^4$ insgesamt keine Quadrate sein können.
\section*{Korollar 5}
Und weil $a^4 - b^4$ kein Quadrat sein kann, werden diese zwei neuen Formeln $\alpha^3 \beta x^4 - \alpha \beta^3 y^4$ und $2\alpha^3\beta x^4 + 2\alpha \beta^3 y^4$ enstehen, die auf keine Weise Quadrate ergeben können.
\section*{Bemerkung}
Aus dem, was wir bis hierher bewiesen haben, werden also die folgenden $6$ allgemeineren Formeln hervorgehen, die auf keine Weise in Quadrate verwandelt werden können:
\[
\begin{array}{llll}
\text{I.} & \alpha^3 \beta x^4 + \alpha \beta^3 y^4 &\qquad \text{IV.} & 2\alpha^3 \beta x^4 - 2\alpha\beta^3 y^4 \\
\text{II.} & \alpha^3 \beta x^4 - \alpha \beta^3 y^4 &\qquad \text{V.} & 2\alpha^3 \beta x^4 + 2\alpha \beta^3 y^4 \\
\text{III.} & \alpha^3 \beta x^4 + 2\alpha \beta^3 y^4 & \qquad \text{VI.} & 2\alpha^3 \beta x^4 - 4\alpha \beta^3 y^4
\end{array}
\]
Und in diesen $6$ Formeln sind alle enthalten, welche wir in den vorhergehenden Formeln behandelt haben. Aus diesen Formeln aber könnten, wie ich es zuvor gemacht habe, trionomiale Formeln gefunden werden, welche in gleicher Weise sicher wären, keinenfalls Quadrate ergeben zu können; aber das Beschaffen dieser erspare ich mir, der ich zu einigen anderen Theoremen fortschreiten will, die sich mit Kuben beschäftigen und die aus diesen Formeln nicht gefunden werden können.
\section*{Theorem 10}
Es kann kein Kubus, nicht einmal nachdem gebrochene Zahlen ausgenommen wurden, um die Einheit vermehrt ein Quadrat ergeben außer dem einzigen Fall, in dem der Kubus $8$ ist.
\section*{Beweis}
Die Proposition geht also darauf zurück, dass $\frac{a^3}{b^3}+1$ niemals ein Quadrat sein kann außer im Fall, in dem $\frac{a}{b} = 2$ ist. Deshalb wird zu zeigen sein, dass diese Form $a^3 b + b^4$ niemals ein Quadrat sein kann, wenn nicht $a = 2b$ ist. Dieser Ausdruck wird aber in diese drei Faktoren $b(a+b)(aa-ab+bb)$ aufgelöst, die zuerst Quadrate festsetzen können, wenn $b(a+b) = aa - ab + bb$ sein könnte, woher $a= 2b$ hervorgeht, welches der Fall ist, den wir ausgenommen haben. Ich setze aber, um weiter fortzufahren, $a+b = c$ oder $a = c-b$, nach Ausführung welcher Substitution man
\[
	bc(cc-3bc+3bb)
\]
haben wird, welcher Ausdruck zu beweisen ist, kein Quadrat sein zu können, wenn nicht $c = 3b$ ist; es sind aber $b$ und $c$ zueinander prime Zahlen. Hier tauchen aber erneut zwei zu betrachtende Fälle auf, je nachdem, ob $c$ ein Vielfaches von $3$ ist oder nicht; in jenem Fall werden die Faktoren $c$ und $cc-3bc+3bb$ den gemeinsamen Teiler $3$ haben, in diesem werden aber alle drei Faktoren zueinander prim sein.

Es sei also zuerst $c$ nicht durch $3$ teilbar; es wird nötig sein, dass jene einzelnen drei Faktoren Quadrate sind, natürlich $b$ und $c$ und $cc-3bc+3bb$ getrennt. Es werde also $cc-3bc+3bb = \left(\frac{m}{n}b - c\right)^2$; es wird
\[
	\frac{b}{c} = \frac{3nn-2mn}{3nn-mm} \quad \text{oder} \quad \frac{b}{c} = \frac{2mn-3nn}{mm-3nn}
\] 
sein, die Terme welches Bruches zueinander prim sein werden, wenn $m$ nicht ein Vielfaches von $3$ ist. Es sei also $m$ durch $3$ nicht teilbar; es wird entweder $c=3nn-mm$ oder $c = mm-3nn$ sein und $b=3nn-2mn$ oder $b=2mn-3nn$. Aber weil $3nn - mm$ kein Quadrat sein kann, setze man $c = mm-3nn$, welches das Quadrat der Wurzel $m - \frac{p}{q}n$ werde, und daher entsteht $\frac{m}{n} = \frac{3qq+pp}{2pq}$ und
\[
	\frac{b}{nn} = \frac{2m}{n} - 3 = \frac{3qq - 3pq + pp}{pq}
\]
Diese Formel $pp(3qq - 3pq + pp)$ wäre also ein Quadrat, welche ganz und gar der vorgelegten Formel $bc(3bb - 3bc + cc)$ gleich ist und aus viel kleineren Zahlen besteht. Aber es sei $m$ ein Vielfaches von $3$, z.\,B. $m=3k$; es wird $\frac{b}{c} = \frac{nn-2kn}{nn-3kk}$ sein, woher entweder $c = nn-3kk$ oder $c = 3kk - nn$ sein wird; weil aber $3kk - nn$ kein Quadrat sein kann, setze man $c = nn-3kk$ und seine Wurzel $n-\frac{p}{q}k$, woher $\frac{n}{k} = \frac{3qq+pp}{2pq}$ oder $\frac{k}{n} = \frac{2pq}{3qq+pp}$ werden wird und
\[
	\frac{b}{nn} = 1 - \frac{2k}{n} = \frac{pp+3qq-4pq}{3qq+pp}.
\]
Es müsste also $(pp+3qq)(p-q)(p-3q)$ ein Quadrat sein. Man setze $p-q = t$ und $p-3q = u$; es wird $q = \frac{t-u}{2}$ und $p = \frac{3t-u}{2}$ sein und jene Formel geht über in diese $tu(3tt-3tu+uu)$, welche wiederum der ersten $bc(3bb-3bc+cc)$ gleich ist. Es bleibt also der letzte Fall übrig, in dem $c$ ein Vielfaches von $3$ ist, z.\,B. $c = 3d$, und es muss $bd(bb-3bd+3dd)$ ein Quadrat sein; weil diese wiederum der ersten gleich ist, ist klar, dass in jedem der beiden Fälle nicht passieren kann, dass die vorgelegte Formel ein Quadrat ist. Deswegen ist außer dem Kubus $8$ nicht einmal in gebrochenen Zahlen ein anderer gegeben, der mit der Einheit ein Quadrat ergibt.

\hfill\textsc{Q.E.D.}
\section*{Korollar 1}
Auf die gleiche Weise kann gezeigt werden, dass kein Kubus um die Einheit vermindert ein Quadrat sein kann; und das nicht einmal in gebrochenen Zahlen.
\section*{Korollar 2}
Daher folgt, dass weder $x^6 + y^6$ noch $x^6 - y^6$ ein Quadrat sein können und keine Dreieckszahl ein Kubus ist außer der Einheit.
\end{document}